\UseRawInputEncoding
\documentclass[12pt]{article}

\usepackage{mathrsfs}
\usepackage{amssymb}
\usepackage{amsmath}
\usepackage{amsfonts}
\usepackage{amstext}
\usepackage{amsthm}
\usepackage{graphicx}
\usepackage{subfig}
\usepackage{color}
\usepackage{indentfirst,latexsym,bm}
\usepackage{mathrsfs}
\usepackage{cite}
\usepackage[all]{xy}
\usepackage{newlfont}

\newcommand{\R}{\mathbb{R}}

\usepackage{enumerate}
\newtheorem{thm}{Theorem}[section]

\numberwithin{equation}{section}

\newtheorem{theorem}{Theorem}[section]

\newtheorem{lemma}[theorem]{Lemma}
\newtheorem{remark}[theorem]{Remark}
\newtheorem{definition}[theorem]{Definition}
\newtheorem{proposition}[theorem]{Proposition}
\numberwithin{equation}{section}

\newcommand{\Sph}{\mathbb{S}^1}
\newcommand{\pd}{\partial}
\newcommand{\dd}{\,\mathrm{d}}
\newcommand{\cG}{\mathcal{G}}
\newcommand{\cU}{\mathcal{U}}
\newcommand{\cR}{\mathcal{R}}

\allowdisplaybreaks
\newcommand{\norm}[1]{\left\Vert#1\right\Vert}

\begin{document}

\title{On the Schiffer and Berenstein conjectures with high-frequency for convex domains in the plane
\thanks{The research of G. Dai is supported by the National Natural Science Foundation of China (Grants No.: 12371110).  The research of J. Wei is partially supported by NSERC of Canada. The research of Y. Zhang is supported by the National Natural Science Foundation of China (Grants No.: 12301133).}}
\author{Guowei Dai
\thanks{School of Mathematical Sciences, Dalian University of Technology, Dalian, 116024, P.R. China. E-mail: daiguowei@dlut.edu.cn} \,,
Yingxin Sun
\thanks{School of Mathematical Sciences, Dalian University of Technology, Dalian, 116024, P.R. China. E-mail: sunyingxin2023@mail.dlut.edu.cn} \,,
Juncheng Wei
\thanks{Corresponding author. Department of Mathematics, The Chinese University of Hong Kong,
Shatin, Hong Kong, China. E-mail: jcwei@math.ubc.ca}  \,,
Yong Zhang
\thanks{School of Mathematical Sciences, Jiangsu University, Zhenjiang, 212013, P.R. China. E-mail: zhangyong@ujs.edu.cn}
}
\date{}
\maketitle

\renewcommand{\abstractname}{Abstract}

\begin{abstract}
In this paper, by introducing two-point stationary-phase amplitude defect, we provide a partial positive answer to the Schiffer and Berenstein conjectures in $\mathbb{R}^2$. More precisely, assuming that a bounded uniformly
convex domain $\Omega \subset \mathbb{R}^2$ has a connected boundary of class $C^{2,\epsilon}$ with $\epsilon \in (0,1)$, we show that if, for some nonzero constant $c_D$, the overdetermined elliptic problem
\begin{equation}
-\Delta u = \alpha u \ \text{ in } \ \Omega, \qquad u = 0 \ \text{ on } \ \partial\Omega, \qquad \frac{\partial u}{\partial \nu} = c_{D} \ \text{ on } \ \partial\Omega  \nonumber
\end{equation}
admits a nontrivial solution corresponding to a large eigenvalue $\alpha$, then the domain $\Omega$ must be a disk. Similarly, we establish that if a domain $\Omega \subset \mathbb{R}^2$ has a connected Lipschitz boundary and the problem
\begin{equation}
-\Delta u = \alpha u \ \text{ in } \ \Omega, \qquad \frac{\partial u}{\partial \nu} = 0 \ \text{ on } \ \partial\Omega, \qquad u = c_{N} \ \text{ on } \ \partial\Omega \nonumber
\end{equation}
has a nontrivial solution corresponding to a large eigenvalue $\alpha$, then $\Omega$ is a disk as well.
\end{abstract}

\emph{Keywords:} Schiffer conjecture; Berenstein conjecture; Overdetermined problem; Stationary phase

\emph{AMS Subjection Classification(2020):} 35N05; 35R35; 43A45; 58J50; 65N25

\tableofcontents

\section{Introduction}
\subsection{Background of Schiffer and Berenstein conjectures}
\bigskip
\quad\,
In the present paper, we investigate the following overdetermined elliptic problems
\begin{equation}\label{overdeterminedeigenvalueproblem}
\left\{
\begin{array}{ll}
\Delta u+\alpha u=0\,\, &\text{in}\,\, \Omega,\\
u=0 &\text{on}\,\, \partial \Omega,\\
\partial_\nu u=c_{D} &\text{on}\,\, \partial \Omega
\end{array}
\right.
\end{equation}
and
\begin{equation}\label{SchifferConjecture}
\left\{
\begin{array}{ll}
\Delta u+\alpha u=0\,\, &\text{in}\,\, \Omega,\\
\partial_\nu u=0 &\text{on}\,\, \partial \Omega,\\
u=c_{N} &\text{on}\,\, \partial \Omega,
\end{array}
\right.
\end{equation}
where $\Omega\subset \mathbb{R}^{N+1}$ for $N\geq 1$ is a bounded domain, $\alpha$ is a parameter, $c_D$ and $c_N$ are two constants and $\nu$ is the unit outer normal vector to $\partial \Omega$.
The study of these overdetermined problems is largely motivated by the works of Schiffer \cite{Schiffer,Schiffer1} and of Berenstein and Yang \cite{Berenstein, Berenstein1, Berenstein2}. Subsequently, increasing attention has been devoted to the following intriguing problem in spectral geometry: Is the ball the unique bounded, smooth, simply connected domain for which problems \eqref{overdeterminedeigenvalueproblem} and \eqref{SchifferConjecture} admit a nontrivial solution? This problem is known as the Schiffer conjecture \cite[Problem 80]{Yau}(or the Berenstein conjecture) and is stated as follows.
\\ \\
\textbf{Schiffer conjecture.} \emph{Let $\Omega\subset \mathbb{R}^{N+1}$ be a bounded regular simply-connected domain. If $u: \Omega \rightarrow \mathbb{R}$
is a nontrivial solution to \eqref{SchifferConjecture},
then $\Omega$ is a ball and $u$ is radially symmetric.}\\
\\
\textbf{Berenstein conjecture.}
\emph{Let $\Omega$ be a bounded $C^{2,\epsilon}$ ($\epsilon\in(0,1)$) domain in $\mathbb{R}^{N+1}$ with connected boundary. If there exists a nontrivial
solution $u$ to the overdetermined eigenvalue problem \eqref{overdeterminedeigenvalueproblem}, then $\Omega$ is a ball.}\\

The overdetermined problem \eqref{SchifferConjecture} is closely related to the Pompeiu problem \cite{Pompeiu}. Specifically, if the only continuous function $f$ satisfying $\int_{\sigma(\Omega)}f(x)\,\text{d}x=0$ for every rigid motion $\sigma$ of $\Omega$ is the zero function, then $\Omega$ is said to have the Pompeiu property. The problem of characterizing which domains possess the Pompeiu property is known as the Pompeiu problem. This problem originated in harmonic analysis and was first proposed by the Polish mathematician Pompeiu in 1929 \cite{Pompeiu1}. The Pompeiu problem has numerous applications in computer-aided tomography \cite{Shepp}, plasma physics \cite{Temam}, nuclear reactor theory \cite{Berenstein, Deans}, satellite remote sensing or geographic mapping \cite{Smith,  Willms}; for comprehensive surveys, see Zalcman \cite{Zalcman} and Chamberland \cite{Chamberland}.
Many domains are known to have the Pompeiu property, including squares \cite{Christov, Pompeiu}, polygons, ellipses \cite{BST}, ellipsoids, certain solid tori in
$\mathbb{R}^N$ \cite{Dalmasso}, any planar domain whose boundary contains a straight line segment or a circular arc (but is not a complete circle) \cite{Williams0}, and convex planar domains whose minimal diameter is at most half of the maximal diameter \cite{BK}, simply-connected bounded planar domains whose boundary can be denoted by some polynomial
with complex coefficients \cite{Garofalo}, quadrature domains of non-disk shapes \cite{Ebenfelt}. Conversely, Tchakaloff (1944) first observed that an open ball of any radius does not have the Pompeiu property \cite{Tchakaloff}. This naturally raises a long-standing open question: must a bounded domain lacking the Pompeiu property necessarily be a ball?
Williams \cite[Theorem 1]{Williams1} proved that for a bounded regular simply-connected domain $\Omega\subset \mathbb{R}^{N+1}$, the domain $\Omega$ fails to have the Pompeiu property if and only if the overdetermined problem (\ref{SchifferConjecture}) admits a nontrivial solution. Consequently, proving the Schiffer conjecture would imply that any bounded regular simply-connected domain without the Pompeiu property must be a ball (in other words, a bounded regular simply connected non-spherical region must possess the Pompeiu property), thereby resolving the Pompeiu problem. The Pompeiu problem or the Schiffer conjecture appears in Yau's list of problems \cite[Problem 80]{Yau} and is also included in Henrot's book \cite[Conjecture 6.57]{Henrot}.

This question may also be related to Kac's problem. In 1966, Kac \cite{Kac} famously asked, ``Can one hear the shape of a drum?" From a mathematical perspective, the goal is to determine the shape of a domain from its eigenvalues, a question to which the Schiffer conjecture directly belongs. The Pompeiu problem is also related to Morera problem because of its relationship to the classical theorem of Morera (see \cite{BST}). Let $\{\Gamma\}$ denote a collection of rectifiable
closed curves in the plane, and let $f$ be continuous on $\mathbb{R}^2$. Suppose
that, setting $\zeta = x + iy$,
\begin{equation*}
\int_{\sigma(\Gamma)}f(\zeta)\,\text{d}\zeta=0,\,\,\Gamma\in\{\Gamma\}.
\end{equation*}
Does this imply that $f$ is an entire holomorphic function of $\zeta$? Zalcman \cite{Zalcman1972} proved that when $\Omega$ is a planar domain whose
boundary $\partial\Omega$ is a rectifiable curve, then an affirmative answer to the
Pompeiu problem for $\Omega$ implies, via Green's Theorem, an affirmative
solution to the Morera problem for $\partial\Omega$. Brown, Schreiber and Taylor \cite{BST} proved that the converse was also true. Therefore, if the Schiffer conjecture holds true, the bounded regular single connected domain without the above Morera property must be a ball.

The problem \eqref{overdeterminedeigenvalueproblem} originates from the shape optimization problem for the principal eigenvalue. Rayleigh \cite{Rayleigh} conjectured that among all fixed membranes of given area, the disk has the smallest principal frequency, which was confirmed by Faber \cite{Faber} and Krahn \cite{Krahn, Krahn1}. Schiffer \cite{Schiffer, Schiffer1} further showed that if the first Dirichlet eigenvalue functional $\Omega\to\lambda_1(\Omega)$ attains a minimum among all bounded smooth domains of prescribed volume, then problem (\ref{overdeterminedeigenvalueproblem}) with $\alpha=\lambda_1$ admits a solution. Serrin's result \cite{Serrin, Pucci} established that the ball is the unique domain for which (\ref{overdeterminedeigenvalueproblem}) admits a positive solution. This result has inspired extensive research \cite{Aftalion, BCN, ChenLi, Del, Gidas, Kamburov, Pucci, Reichel, RRS}, but proving the statement without the positivity assumption is still an open problem \cite{Kawohl}.

\subsection{Previous progress and our main results}

Although these two conjectures remain open in all dimensions, significant partial results have been obtained. For convenience in presenting these results, we denote by $\lambda_k$ and $\mu_k$ the Dirichlet and Neumann eigenvalues, respectively.\\

\textbf{Regarding the Schiffer conjecture:} Berenstein and Yang \cite{Berenstein,Berenstein1,Berenstein2} proved that if the corresponding eigenvalue problem admits infinitely many solutions, then $\Omega$ must be a ball (see \cite{Berenstein} for two-dimensional Euclidean space, \cite{Berenstein1} for the hyperbolic disk, and \cite{Berenstein2} for arbitrary dimensions). Moreover, Berenstein \cite{Berenstein} showed that $\Omega$ is a disk when $\alpha=\lambda_2$ which has been extended into any dimension in \cite{DaiLS}. Aviles \cite{Aviles} extended this to convex domains with $\alpha\leq\mu_7$. These results rely on comparisons between Neumann and Dirichlet eigenvalues. For planar domains, P\'{o}lya and Szeg\H{o} \cite{Polya, Szego} demonstrated that $\mu_2<\lambda_1$ holds in any domain bounded by an analytic curve. Using the idea from \cite{Payne}, Levine and Weinberger \cite{Levine} extended this into $\mathbb{R}^d$ by proving $\mu_{k+r}<\lambda_k$ for any $r\in\{1,\ldots,d\}$, under suitable conditions on the principal curvatures of a
$C^{2+\epsilon}$-smooth boundary of a bounded domain. Deng \cite{Deng} removed the convexity assumption in \cite{Aviles} and further proved, using variational methods, that if $\Omega$ is strictly convex and centrally symmetric with $\alpha<\mu_{13}$, then $\Omega$ is a disk. Willms and Gladwell \cite{WillmsGladwell} proved that $\Omega$ is a ball if it is simply connected and the Neumann eigenfunction $u$ has no interior saddle points. Liu \cite{Liu} established the conjecture under the condition that the third-order interior normal derivative of $u$ is constant on $\partial\Omega$, a result later extended by Kawohl and Marcello \cite{Kawohl} to fourth- or fifth-order normal derivatives. More recently, Mondal \cite{Mondal} proved the Schiffer
conclusion for centrally symmetric planar domains with analytic boundary
when the number of boundary points whose normal lines pass through the
center is less than eight.  This is an all-frequency result under an
analytic and finite-normal-count hypothesis. For further related results, see \cite{Bagchi, Garofalo, Friedlander} and the references therein.

Recently, two surprising counterexamples to problem \eqref{SchifferConjecture} have been obtained by Enciso, Fern\'{a}ndez, Ruiz and Sicbaldi \cite{EncisoFRS} for a bounded planar doubly connected domain, and by Fall, Minlend and Weth \cite{FallMW} for an unbounded cylindrical domain. Two recent preprints announce constructions of counterexamples to the
unrestricted planar Schiffer conjecture.  Colbrook and Stepaniants
\cite{ColbrookStepaniants} gave a computer-assisted real-analytic simply
connected example, while Cao-Labora and De Dios Pont \cite{CaoDeDios} gave
infinitely many smooth $N$-fold symmetric examples. They show why a positive result
should identify a quantitative geometric mechanism that excludes
high-frequency counterexamples in a controlled class.\\

\textbf{Regarding the Berenstein conjecture:} Berenstein \cite{Berenstein} first proved that the existence of infinitely many eigenvalues for problem \eqref{overdeterminedeigenvalueproblem} forces $\Omega\subset\mathbb{R}^2$ to be a disk, which was subsequently extended by Berenstein and Yang to the hyperbolic disk \cite{Berenstein1} and to arbitrary dimensions \cite{Berenstein2}. Liu \cite{Liu} showed that the Berenstein conjecture holds if and only if the second-order interior normal derivative of the corresponding Dirichlet eigenfunction $u$ is constant on $\partial\Omega$. Kawohl and Marcello \cite{Kawohl} generalized this to the case of constant fourth- or fifth-order normal derivatives. Another notable advance concerns the construction of counterexamples on unbounded domains. Following earlier work by Sicbaldi \cite{Sicbaldi}, Dai and Zhang \cite{DaiZ} demonstrated that a perturbation of the cylinder $B_1\times\mathbb{R}$ (where $B_1$ is the unit ball in $\mathbb{R}^n$) yields an unbounded domain with a periodic hypersurface of revolution about the $\mathbb{R}$-axis, on which problem \eqref{overdeterminedeigenvalueproblem} admits a sign-changing solution.\\

In this paper, we provide a partial resolution of both conjectures in the two-dimensional
case.
The following statements are consequences of the main quantitative result Proposition \ref{thm:main-distance}.

\begin{thm}\label{Theorem 1.1}
Let $\Omega$ be a bounded uniformly
convex domain in $\mathbb{R}^{2}$, whose connected boundary $\partial \Omega$ is of class $C^{2,\epsilon}$ with $\epsilon\in(0,1)$. If problem (\ref{overdeterminedeigenvalueproblem}) with $\alpha>\alpha_0(\Omega)$ for some positive constant $\alpha_0(\Omega)$
has an eigenfunction $u$, then $\Omega$ is a disk.
\end{thm}
\begin{thm}\label{Theorem 1.2}
Let $\Omega\subseteq\mathbb{R}^{2}$ be a bounded uniformly
convex domain with a Lipschitz connected boundary $\partial \Omega$. If problem (\ref{SchifferConjecture}) has an eigenfunction $u$ corresponding to a eigenvalue $\alpha>\alpha_0(\Omega)$ for some positive constant $\alpha_0(\Omega)$, then $\Omega$ is a disk.
\end{thm}

From the perspective of physics, the conclusions of Theorems 1.1 and 1.2 imply that one can determine the shape of a drum from a single large tone, provided that the drum satisfies suitable overdetermined boundary conditions. Here, overdetermined boundary conditions mean that both the displacement (height) and its normal derivative are prescribed on the boundary $\partial \Omega$. Under such smooth boundary conditions, a drum for which the boundary values of the function and its normal derivative are known can be uniquely identified as a disk using only one large eigenvalue. This result is currently restricted to two-dimensional case. Its extension to arbitrary dimensions is left for future work.
\subsection{The approach of Berenstein and Yang \cite{Berenstein2}}

The main methods of this paper are rooted in Berenstein and Yang \cite{Berenstein2}. To better present our proof strategy and facilitate comparison with \cite{Berenstein2}, this subsection summarizes the techniques and tools employed in that work. The core goal of their paper \cite{Berenstein2} is to prove that if a bounded smooth domain $\Omega$ (in Euclidean or hyperbolic space) admits a sequence of sufficiently large overdetermined eigenvalues (i.e., corresponding eigenfunctions satisfying problem \eqref{overdeterminedeigenvalueproblem} or \eqref{SchifferConjecture}, then $\Omega$ must be a ball (or a geodesic ball). The main steps of their proof are as follows (taking only the Euclidean space case as an example).
~\\

\textbf{Step 1: Construction of plane wave test functions}

Choose mutually orthogonal unit vectors $\xi$, $\eta$ and define the plane wave solutions
$$
\varphi_{\lambda,t}=e^{i\langle\lambda \xi-i t \eta, x\rangle}=e^{t\eta\cdot x}e^{i\lambda\xi\cdot x},
$$
which satisfies $\Delta\varphi+\left(\lambda^2-t^2\right)\varphi=0$. The parameters $t$ and $\lambda$ are chosen so that the corresponding eigenvalues $\alpha= \lambda^2 - t^2$ tend to $+\infty$, in particular,
$$
|t|\rightarrow+\infty  \quad \text{and} \quad \frac{|t|}{\log\lambda}\rightarrow 0.
$$
\indent
\textbf{Step 2: Green's identity reduces the problem to boundary integrals}

For an overdetermined Dirichlet or Neumann eigenfunction $u$, Green's identity gives
$$
0=\int_{\Omega}\left(u\Delta\varphi-\varphi\Delta u\right)\,\text{d}x=\int_{\partial\Omega}\left(u\frac{\partial\varphi}{\partial\nu}-\varphi\frac{\partial u}{\partial\nu}\right)\,\text{d}S,
$$
$\nu$ is the unit outer normal vector to $\partial \Omega$.
Using the boundary conditions, this leads to the equality of the oscillatory integrals
$$
0=\int_{\partial\Omega}e^{t\eta\cdot x}e^{i\lambda\xi\cdot x}\,\text{d}S \quad \text{or} \quad 0=\int_{\partial\Omega}\langle t \eta+i\lambda \xi, \nu\rangle e^{t\eta\cdot x}e^{i\lambda\xi\cdot x}\,\text{d}S.
$$

\textbf{Step 3: Stationary phase method}

The main contribution to the integral comes from finite number of non-degenerate critical points $p_1, p_2, ...,p_m$ for $\widetilde{\nu}\left(p_{i}\right)=\pm\xi$ with $i\in\{1,\dots,m\}$ of the phase, i.e., points where the unit interior normal $\widetilde{\nu}$ is parallel to $\xi$. Using H\"{o}rmander's asymptotic formula \cite[Theorem 7.7.5]{Hormander} for oscillatory integrals, they obtained that
$$
\sum_{\widetilde{\nu}\left(p_{i}\right)=\pm \xi}\left(2\pi i\lambda\right)^{-\frac{N}{2}}\left(\prod_1^Nk_j(p_{i})\right)^{-\frac{1}{2}}e^{t\eta\cdot p_{i}}e^{i\lambda\xi\cdot p_{i}}+\text{lower order terms},
$$
where $k_j$ are principal curvatures. Since the total integral is zero (from Green's identity), the leading terms will cancel, yielding desired curvature identity.
~\\

\textbf{Step 4: Deriving geometric properties I and II}

The key step is to derive the following two geometric properties:

\textbf{Constant breadth (Property I)}: By varying the parameters $\xi$ and $\eta$ and differentiating with respect to a rotation angle, one obtains a system of linear equations whose coefficients form a Vandermonde matrix. This forces the existence of a pair of points $p$ and $p^*$ such that $\eta\cdot p=\eta\cdot p^*$ for all $\eta$, implying $p^*=p-L\nu$ with constant $L$. Thus $\Omega$ has constant breadth, and the map $p\mapsto p^*$ is an involution.

\textbf{Curvature relation (Property II)}:
Comparing the leading-order coefficients gives $\prod k_j(p)=\pm\prod k_j\left(p^*\right)$.

Using the Property I and Property II, one derives
$/\lambda_{i}+1/\lambda_{i}^*=L$. Combining these yields the curvature identity in Euclidean space
$$
\prod_{i=1}^n \left(L\lambda_{i}-1\right)=\pm 1.
$$

\textbf{Step 5: Proving the domain is a ball}

By analyzing the next term in the asymptotic expansion (involving the $\Box$ operator), one obtains a more refined identity:
$$
\sum_{i=1}^n\frac{1}{k_{i}}\eta_i^2=\sum_{i=1}^n\frac{1}{k_{i}^*}\eta_i^2
$$
for all $\eta$ lying in a small neighborhood in $\eta_0$, which forces $k_{i}(p)=k_{i}\left(p^*\right)$ for each principal curvature. Together with $1/\lambda_{i}+1/\lambda_{i}^*=L$, they deduce that $k_{i}=2/L$ is a constant, so the boundary is a sphere.
~\\

\textbf{Summary}: Transform the inverse problem into an asymptotic analysis of oscillatory integrals. Use the stationary phase method to extract geometric information (constant breadth and curvature relations) from the boundary integrals. Then apply classical differential geometry results (Alexandrov's theorem) to conclude rigidity.

\subsection{Our argumentation ideas and article organization}
The aim of this paper is to prove that, for a two-dimensional bounded uniformly convex domain, the existence of a single large overdetermined eigenvalue forces the domain to be a disk. Our proof is inspired by the work of Berenstein and Yang \cite{Berenstein2} but generalizes it significantly by removing the need for an infinite sequence of eigenvalues $\alpha$ tending to infinity for \eqref{overdeterminedeigenvalueproblem} and \eqref{SchifferConjecture}.
Here we need only a single large eigenvalue $\alpha$ for these problems. Specifically, we express $\alpha$ as $\lambda^2 - t^2$ with $\lambda > 0$ and $t\in \mathbb{R}$. Fixing $t\in\{-1,1\}$, and since $\alpha>0$ is fixed, it follows that $\lambda$ is also a large fixed parameter (rather than tending to infinity).

Firstly, we introduce the two-point amplitude
defect $\mathfrak B(\Omega)$ arising directly from the two stationary
contributions and prove that it vanishes exactly for disks. Second, we compare
this defect quantitatively with the $C^2$ support-function distance on
controlled $C^m$ classes, prove the sharp finite-regularity exponent, and
derive a stronger Hausdorff comparison under central symmetry. Third, we
establish a uniform finite-regularity two-point stationary-phase expansion and
use it to obtain class-uniform high-frequency exclusion for both
overdetermined problems.

We describe the mechanism in more detail. Uniform convexity makes
the phase $x\mapsto\xi\cdot x$ have exactly two nondegenerate stationary
points. The moduli of their leading contributions are
\[
 \sqrt{\rho(\theta)}e^{t h'(\theta)}
 \quad\hbox{and}\quad
 \sqrt{\rho(\theta+\pi)}e^{-t h'(\theta+\pi)}.
\]
Their maximal difference defines $\mathfrak B(\Omega)$. A geometric result
shows that this defect vanishes identically only for a disk. Imposing $\mathfrak B(\Omega)\geq\delta$ makes this
separation uniform over a geometric class. The central quantitative step is
to eliminate this auxiliary assumption from the geometric statement. If
$d_{\mathrm{disk}}(\Omega)$ denotes the $C^2$ support-function distance from the
disks, then for every integer $m\geq5$ we prove on each controlled $C^m$
class that
\[
\mathfrak{B}(\Omega)\geq
 c_m\min\left\{
 d_{\mathrm{disk}}(\Omega)^{(m-1)/(m-2)},1\right\}.
\]
The exponent comes from the even--odd decomposition under the antipodal map,
invertibility of $1+\partial_\theta^2$ on the Steiner-normalized odd
component, and a $C^m$ periodic interpolation inequality for the even
component. Combining this estimate with the two-point stationary-phase
expansion yields an eigenvalue cutoff of order
$1+\varepsilon^{-2(m-1)/(m-2)}$ on families satisfying
$d_{\mathrm{disk}}(\Omega)\geq\varepsilon$. For every finite $m$, a centrally symmetric
high-frequency family shows that the geometric defect--distance exponent
$(m-1)/(m-2)$ is optimal under only uniform $C^m$ control.  This construction
does not assert sharpness of the spectral cutoff.

Central symmetry nevertheless produces a different and genuinely stronger
geometric statement. In Steiner position it makes both $h$ and $\rho$
$\pi$-periodic, so the two-point defect becomes
$2\max\sqrt\rho\,\sinh\left\vert h'\right\vert$. The oscillation of a $\pi$-periodic support
function controls $\left\vert h'\right\vert_\infty$ directly. Consequently, the natural
Hausdorff distance from the centered disks yields a linear lower bound for
the defect and an eigenvalue cutoff of order $1+\varepsilon^{-2}$. This
refinement needs only $C^5$ regularity because no $C^2$ interpolation is
used. It complements, rather than replaces, the nonsymmetric theorem:
central symmetry improves the metric and the cutoff exponent, but it does
not improve the sharp regularity-dependent power when one insists on
measuring $C^2$ distance.

Thus the present results differ from the all-frequency theorem of Mondal:
we remove the finite-normal-count condition. The
conclusion is class-uniform only above a controlled frequency threshold.
Compared with the earlier Fourier-integral
literature, the new point is the two-point defect that characterizes disks
and its quantitative comparison with a natural geometric distance, together
with explicit tracking of the quantities on which the uniform remainder and
threshold depend.

The remainder of this paper is structured as follows. Section 2 introduces some preparatory results.
Specifically, we first define a two-point amplitude defect by means of the support function.
Further, we prove that the defect vanishes if and only if the domain is a disk.
We also introduce the controlled geometric classes.
In Section 3, we complete the proof of Theorem 1.1 using the stationary-phase method and the two-point amplitude defect estimate.
Section 4 presents the proof of Theorem 1.2. Given the structural similarities to Theorem 1.1, we focus primarily on elucidating the essential differences.
The last Section concludes the paper by discussing the sharpness of the regularity-dependent exponent and clarifying the scope of the high-frequency exclusion results.

\section{Preliminary}
\bigskip
\quad\, Prior to the main discussion, this section introduces some preparatory results.
\subsection{Expansion for an oscillatory integral}
We first recall H\"{o}rmander's theorem on oscillatory integrals \cite[Theorem 7.7.5]{Hormander}, which can be stated as follows.
\begin{proposition}[Oscillatory integrals]\label{Proposition 2.1}
Let $U\subset \mathbb{R}^{N}$ be a compact set, $X$ an open neighborhood of $U$ and $k$ a positive integer. If $u \in C_{0}^{2 k}(U), f \in C^{3 k+1}(X)$ and $\operatorname{Im} f \geq 0$ in $X$, $\operatorname{Im} f\left(x_{0}\right)=0, f^{\prime}\left(x_{0}\right)=0, \operatorname{det} f^{\prime \prime}\left(x_{0}\right) \neq 0, f^{\prime} \neq 0$ in $U\backslash\left\{x_{0}\right\}$, then
\begin{equation*}
  \begin{array}{l}
\left|\int_U u(x) e^{i\lambda f(x)} \,\text{d} x-e^{i\lambda f\left(x_{0}\right)}\left(\operatorname{det}\left(\lambda f^{\prime \prime}\left(x_{0}\right) / 2 \pi i\right)\right)^{-\frac{1}{2}}\sum\limits_{j<k} \lambda^{-j} L_{j} u\right| \leq C\lambda^{-k} \sum\limits_{|\alpha| \leq 2 k} \sup \left|D^{\alpha} u\right|
\end{array}
\end{equation*}
\emph{for $\lambda>0$.
Here $C$ is bounded when $f$ stays in a bounded set in $C^{3 k+1}(X)$ and $\left|x-x_{0}\right| /\left|f^{\prime}(x)\right|$ has a uniform bound. With}
\begin{equation*}
  g_{x_{0}}(x)=f(x)-f\left(x_{0}\right)-\left\langle f^{\prime \prime}\left(x_{0}\right)\left(x-x_{0}\right), x-x_{0}\right\rangle / 2
\end{equation*}
\emph{which vanishes of third order at $x_{0}$ we have}
\begin{equation*}
L_{j} u=\sum_{\nu-\mu=j} \sum_{2\nu \geq 3 \mu} i^{-j} 2^{-\nu}\left\langle f^{\prime \prime}\left(x_{0}\right)^{-1} D, D\right\rangle^{\nu}\left(g_{x_{0}}^{\mu} u\right)\left(x_{0}\right) / \mu!\nu!,
\end{equation*}
where $i$ is the imaginary unit. This is a differential operator of order $2j$ acting on $u$ at $x_{0}$.
The coefficients are rational homogeneous functions of degree $-j$ in $f^{\prime \prime}\left(x_{0}\right), \ldots, f^{(2 j+2)}\left(x_{0}\right)$ with denominator $\left(\operatorname{det} f^{\prime \prime}\left(x_{0}\right)\right)^{3 j}$.
In each term the total number of derivatives of $u$ and of $f^{\prime \prime}$ is at most $2j$.
\end{proposition}

The symbol $\operatorname{Im}$ in Proposition~\ref{Proposition 2.1} represents the imaginary part, and $\alpha$ is a multi-index (not to be confused with the eigenvalue). Let $f$ be a real-valued $C^{\infty}$ function defined in a compact neighborhood $U$ of the origin in $\mathbb{R}^{N}$, assuming $f$ has only a non-degenerate critical point at the origin in $U$. By applying an orthogonal transformation and using the Taylor's expansion of the function $f$, we have
\begin{equation*}
f(x)=f(0)+\sum_{i=1}^{N} \frac{1}{2} k_{i} x_{i}^{2}+\sum_{i,j,k=1}^{N} A_{i j k} x_{i} x_{j} x_{k}+O\left(|x|^{4}\right),
\end{equation*}
where $k_{i}\neq0$ with $i\in\{1,\ldots,N\}$ are the principal curvatures of $U$ in $x_{i}$ direction at zero, $A_{i j k}$ is symmetric in all three indices.
Write
\begin{equation*}
g(x)=f(x)-f(0)-\sum_{i=1}^{N} \frac{k_{i}}{2} x_{i}^{2}=\sum_{i,j,k=1}^{N} A_{i j k} x_{i} x_{j} x_{k}+O\left(|x|^{4}\right).
\end{equation*}
For $u\in C_{0}^{\infty}(U)$, applying Proposition~\ref{Proposition 2.1}, Berenstein and Yang \cite{Berenstein2} obtained that
\begin{equation}\label{Hormander}
\begin{array}{c}
\left|\int_U u(x) e^{i \lambda f(x)}\, \text{d}x-e^{i \lambda f(0)}\left(\prod\limits_{j=1}^{N} \frac{\lambda k_{j}}{2 \pi i}\right)^{-\frac{1}{2}}\left(u(0)+\frac{L(u)}{\lambda}\right)\right|
\leq C \lambda^{-\frac{N}{2}-2} \sum\limits_{|\alpha| \leq 4} \sup \left|D^{\alpha} u\right|
\end{array}
\end{equation}
for $\lambda>0$, where $D$ represents the first-order derivative operator and $C$ is bounded when $f$ within a bounded set in $C^{7}(U)$ and $|x|/\left|f^{\prime}(x)\right|$ is bounded,
\begin{equation*}
L(u)=-i\left[\frac{1}{2} \Box u(0)+\frac{1}{2^{2}}\frac{1}{ 2 !} \Box^{2}(g u)(0)+\frac{1}{2^{3}} \frac{1}{2 !} \frac{1}{3 !} \Box^{3}\left(g^{2} u\right)(0)\right]
\end{equation*}
with
\begin{equation*}
\Box=\sum_{i=1}^{N} \frac{1}{k_{i}} \frac{\partial^{2}}{\partial x_{i}^{2}}.
\end{equation*}
If $f$ has a finite number of non-degenerate critical points, a suitable partition of unity can reduce the integral to the same finite sum as the integral of the above type when applying this expansion over a compact hyper-surface. Outside the critical points, the integral decreases faster than $\lambda^{-k}$ for any $k$.

Formula (\ref{Hormander}) is a crucial tool used by Berenstein and Yang \cite{Berenstein2}. In their work, the eigenvalue
$\alpha$ is represented as $\lambda^2-t^2$. To apply (\ref{Hormander}) and derive their desired conclusions, they assumed that
\begin{equation}\label{Berensteincondition}
|t|\rightarrow+\infty\quad \text{and} \quad \frac{|t|}{\log \lambda} \to 0.
\end{equation}
Without this assumption, their argument would collapse. However, this assumption (\ref{Berensteincondition}) holds only when
$\alpha\rightarrow+\infty$, which inherently presupposes the existence of a sequence of eigenvalues.
In our case, where only a single eigenvalue is assumed to exist, the assumption (\ref{Berensteincondition}) fails. To overcome this difficulty, we need to find new methods.
\subsection{Support functions and the geometric class}

We use the Steiner point to fix translations canonically. If $h$ is the
support function of a convex body, its Steiner point in the plane is
\begin{equation}\label{eq:Steiner}
s(\Omega)=\frac1\pi\int_0^{2\pi}h(\theta)n(\theta)\,d\theta.
\end{equation}
After translating $s(\Omega)$ to the origin ($\Omega$ is called in Steiner position), set
\begin{equation*}
n(\theta)=(\cos\theta,\sin\theta),\qquad
\tau(\theta)=(-\sin\theta,\cos\theta).
\end{equation*}
For a uniformly convex domain, define the support function
\begin{equation*}\label{eq:support}
h(\theta)=\max_{x\in\overline\Omega}x\cdot n(\theta).
\end{equation*}
The boundary point with outer normal $n(\theta)$ is
\begin{equation}\label{eq:p}
p(\theta)=h(\theta)n(\theta)+h'(\theta)\tau(\theta).
\end{equation}
Its radius of curvature and curvature are
\begin{equation*}\label{eq:rho}
\rho(\theta)=h(\theta)+h''(\theta)>0,\qquad
\kappa(\theta)=\rho(\theta)^{-1}.
\end{equation*}

Based on the support function, we introduce the concept of two-point amplitude defect.
\begin{definition}[Two-point amplitude defect]\label{def:defect}
\emph{For a uniformly convex domain in Steiner position, define}
\begin{equation}\label{eq:defect}
\mathfrak B(\Omega)
:=
\max_{\substack{\theta\in\R\\ t\in\{-1,1\}}}
\left|
\sqrt{\rho(\theta)}e^{t h'(\theta)}
-\sqrt{\rho(\theta+\pi)}e^{-t h'(\theta+\pi)}
\right|.
\end{equation}
\end{definition}

The two quantities in \eqref{eq:defect} are exactly the moduli of the two
leading stationary-phase contributions. The following observation shows
that the defect is a genuine quantitative measure of noncircularity rather
than an auxiliary analytic assumption.

\begin{proposition}[The defect characterizes disks]\label{prop:defect-disk}
Let $\Omega$ be uniformly convex and let its support function be of class
$C^2$. Then
\[
 \mathfrak B(\Omega)=0
 \quad\Longleftrightarrow\quad
 \Omega\text{ is a disk}.
\]
\end{proposition}

\begin{proof}
A disk in Steiner position has constant support function, so its defect
vanishes. Conversely, suppose that $\mathfrak B(\Omega)=0$. Applying
\eqref{eq:defect} with $t=1$ and $t=-1$ gives, for every $\theta$,
\begin{align*}
\sqrt{\rho(\theta)}e^{h'(\theta)}
&=\sqrt{\rho(\theta+\pi)}e^{-h'(\theta+\pi)},\\
\sqrt{\rho(\theta)}e^{-h'(\theta)}
&=\sqrt{\rho(\theta+\pi)}e^{h'(\theta+\pi)}.
\end{align*}
Dividing and multiplying these identities yields
\begin{equation}\label{eq:defect-zero-identities}
 h'(\theta+\pi)=-h'(\theta),
 \qquad \rho(\theta+\pi)=\rho(\theta).
\end{equation}
The first identity implies that
$h(\theta)+h(\theta+\pi)=w$ for some constant $w$. Therefore
$h''(\theta+\pi)=-h''(\theta)$, and the second identity in
\eqref{eq:defect-zero-identities} gives
\[
 w-h(\theta)-h''(\theta)=h(\theta)+h''(\theta).
\]
Hence $h+h''=w/2$, whose $2\pi$-periodic solutions are
\[
 h(\theta)=\frac w2+a\cos\theta+b\sin\theta.
\]
The last two terms represent a translation. Since the Steiner point is the
origin, \eqref{eq:Steiner} gives $a=b=0$. Thus $h$ is constant and $\Omega$
is a disk.
\end{proof}

To apply Proposition \ref{prop:defect-disk}, we introduce the following controlled geometric class.
\begin{definition}[Controlled geometric class]\label{def:class}
\emph{Given an integer $m\geq5$ and $R,M,\kappa_0,\kappa_1>0$, let
$\cG_m\left(R,M,\kappa_0,\kappa_1\right)$ be the class of domains satisfying
\begin{enumerate}[(i)]
\item $\Omega$ is in Steiner position and its support function belongs to
$C^m\left(\R/2\pi\mathbb Z\right)$ with $\rho=h+h''>0$;
\item $\Omega\subset B_R(0)$ and
$\norm{h}_{C^m\left(\R/2\pi\mathbb Z\right)}\le M$;
\item $0<\kappa_0\le\kappa(\theta)\le\kappa_1$ for every $\theta$.
\end{enumerate}}
\end{definition}

For $\delta>0$, set
\[
 \cU_m\left(R,M,\kappa_0,\kappa_1,\delta\right)
 :=\left\{\Omega\in\cG_m\left(R,M,\kappa_0,\kappa_1\right):
 \mathfrak B(\Omega)\geq\delta\right\}.
\]
The support-function assumptions imply that $\partial\Omega$ is at least
$C^{m-1}$. Membership in $\cU_m$ additionally imposes a quantitative
noncircularity condition, by Proposition \ref{prop:defect-disk}.  The
curvature bounds imply
\begin{equation}\label{eq:rho-bounds}
\frac1{\kappa_1}\le\rho(\theta)\le\frac1{\kappa_0}.
\end{equation}
Throughout, we use the nested norm
\[
 \|v\|_{C^k}:=\max_{0\le j\le k}\left\|\partial_\theta^jv\right\|_{L^\infty};
\]
replacing it by an equivalent standard $C^k$ norm changes only the displayed
constants.  With this convention, $\cG_m\left(R,M,\kappa_0,\kappa_1\right)$ is
contained in $\cG_5\left(R,M,\kappa_0,\kappa_1\right)$ for every $m\geq5$. Thus the
analytic estimates proved below at the $C^5$ level apply uniformly to every
class in the regularity scale. The purely geometric defect estimates depend
on $M,\kappa_0,$ and $\kappa_1$, but not on $R$. In the analytic module we
retain $R$ in the constants because $\Omega\subset B_R$ gives the direct bound
$|p|\le R$ for the exponential amplitudes. Since \eqref{eq:p} also bounds
$p$ in terms of the nested $C^m$ norm, $R$ could instead be absorbed into
$M$; keeping it separate makes the geometric and analytic dependencies
transparent.

\subsection{Exact boundary-transform identities}

Let $\xi,\eta\in\Sph$ be mutually orthogonal and, for $t\in\R$ and
$\lambda>0$, define the complex plane wave
\begin{equation*}\label{eq:plane-wave}
\phi_{\lambda,t}(x)
=\exp\left(t\eta\cdot x+i\lambda\xi\cdot x\right).
\end{equation*}
Since $\xi\cdot\eta=0$,
\[
\Delta\phi_{\lambda,t}
=\left(t^2-\lambda^2\right)\phi_{\lambda,t}.
\]
Thus $\phi_{\lambda,t}$ solves the same Helmholtz equation as $u$ precisely
when
\begin{equation}\label{eq:dispersion}
\lambda^2-t^2=\alpha.
\end{equation}

\begin{lemma}[Green identities]\label{prop:green}
Assume \eqref{eq:dispersion}.
\begin{enumerate}[(i)]
\item If \eqref{overdeterminedeigenvalueproblem} has a nontrivial solution, then
\begin{equation*}\label{eq:transform-D}
\int_{\pd\Omega}\phi_{\lambda,t}\dd s=0.
\end{equation*}
\item If \eqref{SchifferConjecture} has a nontrivial solution, then
\begin{equation*}\label{eq:transform-N}
\int_{\pd\Omega}\pd_\nu\phi_{\lambda,t}\dd s=0.
\end{equation*}
\end{enumerate}
For a fixed eigenvalue $\alpha$, both identities are valid for every
$t\in\R$ with $\lambda=\left(\alpha+t^2\right)^{1/2}$.
\end{lemma}

\begin{proof}
Green's identity gives
\begin{align*}
0
&=\int_\Omega
\bigl(u\Delta\phi_{\lambda,t}
-\phi_{\lambda,t}\Delta u\bigr)\dd x\\
&=\int_{\pd\Omega}
\bigl(u\pd_\nu\phi_{\lambda,t}
-\phi_{\lambda,t}\pd_\nu u\bigr)\dd s.
\end{align*}
For \eqref{overdeterminedeigenvalueproblem}, the boundary integral equals
\[
-c_D\int_{\pd\Omega}\phi_{\lambda,t}\dd s.
\]
For \eqref{SchifferConjecture}, it equals
\[
c_N\int_{\pd\Omega}\pd_\nu\phi_{\lambda,t}\dd s.
\]
Cauchy
uniqueness implies $c_D\ne0$ and $c_N\ne0$, which proves the identities.
\end{proof}

\begin{remark}
\emph{The relation \eqref{eq:dispersion} is maintained throughout. In the proofs
below we choose $t\in\{-1,1\}$ once and set
\[
\lambda=\sqrt{\alpha+1}.
\]
Thus $t$ remains bounded and fixed, while
$\lambda/\sqrt\alpha\to1$ as $\alpha\to\infty$.  There is no independent
variation of $t$ and $\lambda$, and hence no hidden regime in which the
exponential parameter and the oscillatory parameter tend to infinity
together.}
\end{remark}

\section{Proof of Theorem 1.1}

\subsection{High-frequency exclusion}

For a domain in Steiner position, define its support-function distance from
the family of disks by
\begin{equation}\label{eq:disk-distance}
 d_{\mathrm{disk}}(\Omega)
 :=\inf_{r>0}\|h-r\|_{C^2(\R/2\pi\mathbb Z)}.
\end{equation}
Both $\Omega$ and the comparison disks in \eqref{eq:disk-distance} are
Steiner-normalized.  Thus \eqref{eq:disk-distance} is the canonical
support-function distance from disks modulo translations; the missing first
sine and cosine modes are precisely the removed translation directions.
For an integer $m\geq5$, write
\begin{equation*}\label{eq:p-m}
 p_m:=\frac{m-1}{m-2}.
\end{equation*}

\begin{lemma}[Periodic interpolation]\label{lem:periodic-interpolation}
For every integer $m\geq3$ there exists $C_m>0$ such that every
$u\in C^m\left(\R/2\pi\mathbb Z\right)$ satisfies
\begin{equation}\label{eq:periodic-interpolation}
 \left\|u''\right\|_{L^\infty}
 \leq C_m\left(
 \left\|u'\right\|_{L^\infty}^{(m-2)/(m-1)}
 \left\|u^{(m)}\right\|_{L^\infty}^{1/(m-1)}
 +\left\|u'\right\|_{L^\infty}\right).
\end{equation}
\end{lemma}

\begin{proof}
Apply the one-dimensional Gagliardo--Nirenberg inequality to $g=u'$ on a
fixed finite covering of the circle by coordinate intervals:
\[
 \left\|g'\right\|_{L^\infty}
 \leq C_m\left(
 \left\|g\right\|_{L^\infty}^{1-1/(m-1)}
 \left\|g^{(m-1)}\right\|_{L^\infty}^{1/(m-1)}
 +\left\|g\right\|_{L^\infty}\right).
\]
The covering and its cutoffs are fixed, so the constant depends only on
$m$. Substituting $g=u'$ gives \eqref{eq:periodic-interpolation}.
\end{proof}

Using the support-function distance, we can establish the following defect estimate.
\begin{proposition}[Regularity-dependent defect estimate]
\label{prop:defect-distance}
Let $m\geq5$ be an integer. For every $M,\kappa_0,\kappa_1>0$ there exists
$c_m=c_m\left(M,\kappa_0,\kappa_1\right)>0$ such that, for every $R>0$ and every
$\Omega\in\cG_m\left(R,M,\kappa_0,\kappa_1\right)$,
\begin{equation}\label{eq:defect-distance-lower}
 \mathfrak B(\Omega)
 \geq c_m
 \min\left\{d_{\mathrm{disk}}(\Omega)^{p_m},1\right\}.
\end{equation}
The constant is independent of $R$ and of the particular domain.
\end{proposition}

\begin{proof}
Write $\beta=\mathfrak B(\Omega)$ and abbreviate
$\rho_\pi(\theta)=\rho(\theta+\pi)$ and
$h'_\pi(\theta)=h'(\theta+\pi)$. The four quantities
\[
 X_\pm=\frac12\log\rho(\theta)\pm h'(\theta),
 \qquad
 Y_\pm=\frac12\log\rho_\pi(\theta)\mp h'_\pi(\theta)
\]
range in a compact interval determined only by
$M,\kappa_0,$ and $\kappa_1$. By the definition of $\beta$,
\begin{equation}\label{zhishubound}
 \left|e^{X_\pm}-e^{Y_\pm}\right|\leq\beta.
\end{equation}
Here and below, the constants $C_j$ depend only on
$m,M,\kappa_0,$ and $\kappa_1$.
All four exponentials are bounded below by
$\mu=\kappa_1^{-1/2}e^{-M}>0$.
The logarithmic function $x \mapsto \log x$ is Lipschitz on $[\mu, \infty)$ with constant $\mu^{-1}$.
Therefore, from (\ref{zhishubound}),
\begin{align*}
\left|X_+ - Y_+\right|
&= \left| \log e^{X_+} - \log e^{Y_+} \right|
\le \frac{1}{\mu} \left|e^{X_+} - e^{Y_+}\right|
\le \frac{\beta}{\mu}, \\
\left|X_- - Y_-\right|
&= \left| \log e^{X_-} - \log e^{Y_-} \right|
\le \frac{1}{\mu} \left|e^{X_-} - e^{Y_-}\right|
\le \frac{\beta}{\mu}.
\end{align*}
Adding and subtracting the two resulting inequalities give
\[
\left| \log \rho(\theta) - \log \rho(\theta+\pi) \right| \le \frac{2\beta}{\mu}.
\]
and
\[
\left| h'(\theta) + h'(\theta+\pi) \right| \le \frac{\beta}{\mu}.
\]
Since $\theta$ is arbitrary, taking the $L^\infty$ norm over $\theta \in \mathbb{R}/2\pi\mathbb{Z}$ gives
\[
\left\| h'(\cdot) + h'(\cdot+\pi)\right\|_{L^\infty} \le \frac{1}{\mu}\beta
\]
and
\[
\left\| \log \rho(\cdot) - \log \rho(\cdot+\pi)\right\|_{L^\infty} \le \frac{2}{\mu}\beta.
\]
Adding the two inequalities yields
\begin{equation}\label{eq:log-defect-control}
 \left\|h'(\cdot)+h'(\cdot+\pi)\right\|_{L^\infty}
 +\left\|\log\rho(\cdot)-\log\rho(\cdot+\pi)\right\|_{L^\infty}
 \leq C_0\beta.
\end{equation}

From the curvature bounds in \eqref{eq:rho-bounds}, we have that
both $\rho(\theta)$ and $\rho(\theta+\pi)$ lie in the compact interval
$\left[1/\kappa_1, 1/\kappa_0\right]$ for any fixed $\theta$.
By the mean value theorem applied to $\log x$ on the interval $\left[1/\kappa_1, 1/\kappa_0\right]$,
there exists $\xi$ between $\rho(\theta)$ and $\rho(\theta+\pi)$ such that
\[
\log \rho(\theta) - \log \rho(\theta+\pi) = \frac{1}{\xi}(\rho(\theta)-\rho(\theta+\pi)).
\]
Taking absolute values gives
\[
\left|\rho(\theta)-\rho(\theta+\pi)\right| = \xi |\log \rho(\theta) - \log \rho(\theta+\pi)|.
\]
From (\ref{eq:log-defect-control}) and the fact that $\xi \le 1/\kappa_0$,
we have
\begin{equation*}
|\rho(\theta)-\rho(\theta+\pi)| \le \frac{1}{\kappa_0} |\log \rho(\theta) - \log\rho(\theta+\pi)|\le \frac{C_0}{\kappa_0} \beta.
\end{equation*}
Since the above inequality holds for every $\theta$, taking the supremum over
$\theta \in \mathbb{R}/2\pi\mathbb{Z}$ yields
\begin{equation}\label{eq:rho-defect-control}
 \|\rho(\cdot)-\rho(\cdot+\pi)\|_{L^\infty}
 \leq C_1\beta.
\end{equation}

Decompose $h$ into its antipodally even and odd parts,
\[
 h_{\mathrm{e}}(\theta)=\frac{h(\theta)+h(\theta+\pi)}{2},
 \qquad
 h_{\mathrm{o}}(\theta)=\frac{h(\theta)-h(\theta+\pi)}{2}.
\]
Equations \eqref{eq:log-defect-control}--\eqref{eq:rho-defect-control}
give
\begin{equation}\label{eq:even-odd-control}
 \left\|h_{\mathrm{e}}'\right\|_{L^\infty}\leq C_2\beta,
 \qquad
 \left\|h_{\mathrm{o}}+h_{\mathrm{o}}''\right\|_{L^\infty}\leq C_2\beta.
\end{equation}
Steiner normalization means that the first sine and cosine Fourier
coefficients of $h$, and hence of $h_{\mathrm{o}}$, vanish. The kernel of
$1+\partial_\theta^2$ on $2\pi$-periodic functions is precisely the span of
$\cos\theta$ and $\sin\theta$. The elementary periodic Green-function
estimate on the complementary subspace therefore yields
\begin{equation}\label{eq:odd-C2-control}
 \left\|h_{\mathrm{o}}\right\|_{C^2}\leq
 C_3\left\|h_{\mathrm{o}}+h_{\mathrm{o}}''\right\|_{L^\infty}
 \leq C_4\beta.
\end{equation}

Let $r=(2\pi)^{-1}\int_0^{2\pi}h_{\mathrm{e}}(\theta)\,\mathrm{d}\theta$; this number is
positive because the Steiner point lies in the interior of $\Omega$.
The periodic Poincar\'e inequality and \eqref{eq:even-odd-control} give
\begin{equation*}\label{eq:even-C1-control}
 \left\|h_{\mathrm{e}}-r\right\|_{C^1}\leq C_5\beta.
\end{equation*}
Applying Lemma \ref{lem:periodic-interpolation} to $h_{\mathrm{e}}$ gives
\begin{align}
 \left\|h_{\mathrm{e}}''\right\|_{L^\infty}
 &\leq C_6\left(
 \left\|h_{\mathrm{e}}'\right\|_{L^\infty}^{(m-2)/(m-1)}
 \left\|h_{\mathrm{e}}^{(m)}\right\|_{L^\infty}^{1/(m-1)}
 +\left\|h_{\mathrm{e}}'\right\|_{L^\infty}\right)\nonumber\\
 &\leq C_7\left(\beta^{1/p_m}+\beta\right).
 \label{eq:even-C2-interpolation}
\end{align}
Combining \eqref{eq:odd-C2-control}--\eqref{eq:even-C2-interpolation} and $h=h_{\mathrm{e}}+h_{\mathrm{o}}$ gives
\begin{equation*}\label{eq:distance-upper-by-defect}
 d_{\mathrm{disk}}(\Omega)
 \leq\|h-r\|_{C^2}
 \leq C_8\left(\beta^{1/p_m}+\beta\right).
\end{equation*}
If $\beta\leq1$, this implies
$\beta\geq\left(2C_8\right)^{-p_m}d_{\mathrm{disk}}(\Omega)^{p_m}$; if
$\beta>1$, the required estimate is immediate after decreasing the
constant. This proves \eqref{eq:defect-distance-lower}.
\end{proof}

\begin{proposition}[Regularity-dependent high-frequency exclusion]
\label{thm:main-distance}
Let $m\geq5$ be an integer. For every $R,M,\kappa_0,\kappa_1>0$ there
exists $C_m^{\mathrm{hf}}=C_m^{\mathrm{hf}}\left(R,M,\kappa_0,\kappa_1\right)>0$ such that, for
every $\varepsilon>0$, no domain
$\Omega\in\cG_m\left(R,M,\kappa_0,\kappa_1\right)$ with
\[
 d_{\mathrm{disk}}(\Omega)\geq\varepsilon
\]
admits a nontrivial solution of either \eqref{overdeterminedeigenvalueproblem} or \eqref{SchifferConjecture} when
\begin{equation}\label{eq:epsilon-threshold}
 \alpha>C_m^{\mathrm{hf}}\left(1+\varepsilon^{-2p_m}\right).
\end{equation}
\end{proposition}

The following defect-based result is the analytic core of the paper. For the sake of fluency in writing, we have postponed its proof to the last subsection of this section. The
distance formulations are obtained from it by the geometric conversions in
Propositions \ref{prop:defect-distance} and \ref{prop:symmetric-defect}.

\begin{proposition}[Defect-based uniform exclusion]\label{thm:uniform}
Let $m\geq5$ be an integer. For every
$R,M,\kappa_0,\kappa_1,\delta>0$ there exists
\[
A=A\left(R,M,\kappa_0,\kappa_1,\delta\right)>0
\]
such that no domain
\[
\Omega\in\cU_m\left(R,M,\kappa_0,\kappa_1,\delta\right)
\]
admits a nontrivial solution of either \eqref{overdeterminedeigenvalueproblem} or \eqref{SchifferConjecture} with
$\alpha>A$.
\end{proposition}

\begin{proof}[Proof of Proposition \ref{thm:main-distance} from Proposition
\ref{thm:uniform}]
Proposition \ref{prop:defect-distance} gives
\[
 \mathfrak B(\Omega)\geq
 \delta_\varepsilon:=c_m\min\{\varepsilon^{p_m},1\}.
\]
Apply Proposition \ref{thm:uniform} with the same $m$ and
$\delta=\delta_\varepsilon$.  The
explicit formula \eqref{eq:A-uniform} below gives
\[
 A\left(R,M,\kappa_0,\kappa_1,\delta_\varepsilon\right)
 \leq C_m^{\mathrm{hf}}\left(1+\varepsilon^{-2p_m}\right),
\]
after increasing a constant depending only on $m$ and the four geometric
controls.
This is \eqref{eq:epsilon-threshold}.
\end{proof}

The preceding proposition is formulated without symmetry and uses a $C^2$
distance. We next isolate the improvement that follows from central
symmetry. After translating the center of symmetry to the origin, the
Steiner point is automatically the origin and $h$ is $\pi$-periodic. The
appropriate weaker distance is then the Hausdorff distance from the
\emph{centered} disks.

\begin{definition}[Centered-disk Hausdorff distance]\label{def:H-distance}
\emph{For a centrally symmetric domain in Steiner position, set
\begin{equation*}\label{eq:H-distance}
 d_{H,{\mathrm{disk}}}(\Omega)
 :=\inf_{r>0}\|h-r\|_{L^\infty(\R/2\pi\mathbb Z)}.
\end{equation*}
Equivalently, this is the Hausdorff distance from $\Omega$ to the disks
centered at the origin.}
\end{definition}

\begin{proposition}[Symmetric defect identity and Hausdorff control]
\label{prop:symmetric-defect}
Let $\Omega=-\Omega$ be uniformly convex and in Steiner position. Then
\begin{equation}\label{eq:symmetric-defect-identity}
 \mathfrak B(\Omega)
 =2\max_{\theta}\sqrt{\rho(\theta)}\,\sinh|h'(\theta)|.
\end{equation}
If $\kappa\leq\kappa_1$, then
\begin{equation}\label{eq:defect-symmetric-lower}
 \mathfrak B(\Omega)
 \geq\frac{2}{\sqrt{\kappa_1}}
 \sinh\!\left(\frac{4}{\pi}d_{H,{\mathrm{disk}}}(\Omega)\right).
\end{equation}
In particular, the defect controls the Hausdorff distance linearly near the
family of centered disks.
\end{proposition}

\begin{proof}
Central symmetry gives
$h(\theta+\pi)=h(\theta)$ and
$\rho(\theta+\pi)=\rho(\theta)$. Substitution in
\eqref{eq:defect} proves \eqref{eq:symmetric-defect-identity}.
Since $h$ is positive and $\pi$-periodic,
\begin{equation*}\label{eq:H-oscillation}
 d_{H,{\mathrm{disk}}}(\Omega)
 =\frac12\bigl(\max h-\min h\bigr).
\end{equation*}
On the circle $\R/\pi\mathbb Z$, points at which the maximum and minimum
are attained can be joined by an arc of length at most $\pi/2$. The mean
value theorem therefore gives
\begin{equation*}\label{eq:H-derivative}
 \left\|h'\right\|_{L^\infty}
 \geq\frac{2}{\pi}\left(\max h-\min h\right)
 =\frac{4}{\pi}d_{H,{\mathrm{disk}}}(\Omega).
\end{equation*}
Combining this inequality, \eqref{eq:symmetric-defect-identity}, and
$\rho\geq\kappa_1^{-1}$ proves \eqref{eq:defect-symmetric-lower}.
\end{proof}

\begin{proposition}[Centrally symmetric Hausdorff refinement]
\label{thm:symmetric-H}
For every $R,M,\kappa_0,\kappa_1>0$ there exists
$C_{\mathrm{sym}}=C_{\mathrm{sym}}\left(R,M,\kappa_0,\kappa_1\right)>0$ with the following
property. For every $\varepsilon>0$, no centrally symmetric domain
$\Omega\in\cG_5\left(R,M,\kappa_0,\kappa_1\right)$ satisfying
\[
 d_{H,{\mathrm{disk}}}(\Omega)\geq\varepsilon
\]
admits a nontrivial solution of either \eqref{overdeterminedeigenvalueproblem} or \eqref{SchifferConjecture} when
\begin{equation*}\label{eq:symmetric-polynomial-threshold}
 \alpha>C_{\mathrm{sym}}\left(1+\varepsilon^{-2}\right).
\end{equation*}
More precisely, exclusion holds as soon as
\begin{equation}\label{eq:symmetric-exact-threshold}
 \alpha>
 \left[
 \max\left\{3,
 \frac{C_*\left(R,M,\kappa_0,\kappa_1\right)\sqrt{\kappa_1}}
 {2\sqrt{2\pi}\,\sinh\left(4\varepsilon/\pi\right)}
 \right\}
 \right]^2-1.
\end{equation}
\end{proposition}

\begin{proof}
Proposition \ref{prop:symmetric-defect} gives
\[
 \mathfrak B(\Omega)\geq
 \delta_{\mathrm{sym}}(\varepsilon):=
 \frac{2}{\sqrt{\kappa_1}}\sinh\left(4\varepsilon/\pi\right).
\]
Apply Proposition \ref{thm:uniform} with $m=5$ and
$\delta=\delta_{\mathrm{sym}}(\varepsilon)$. Formula \eqref{eq:A-uniform}
is exactly \eqref{eq:symmetric-exact-threshold}. Finally,
$\sinh x\geq x$ for $x\geq0$, so the right-hand side of
\eqref{eq:symmetric-exact-threshold} is bounded above by
$C_{\mathrm{sym}}\left(1+\varepsilon^{-2}\right)$.
\end{proof}

The Hausdorff improvement does not mean that central symmetry improves the
$C^2$ exponent in Proposition \ref{prop:defect-distance}. The next family
separates these two assertions and records the optimality of that exponent.

\begin{proposition}[Sharpness at every finite regularity]
\label{prop:scale-sharpness}
Let $m\geq5$ be an integer and fix a sufficiently small $\eta>0$. For even
integers $N\geq2$, let
\begin{equation*}\label{eq:sharp-family}
 h_{m,N}(\theta)=1+\eta N^{-m}\cos(N\theta).
\end{equation*}
Then $h_{m,N}$ are support functions of centrally symmetric uniformly convex
domains $\Omega_{m,N}$ belonging to one controlled $C^m$ class, and
\begin{equation}\label{eq:sharp-asymptotics}
 d_{\mathrm{disk}}\left(\Omega_{m,N}\right)\asymp_\eta N^{-(m-2)},
 \qquad
 \mathfrak B\left(\Omega_{m,N}\right)\asymp_\eta N^{-(m-1)}
 \asymp_{m,\eta} d_{\mathrm{disk}}\left(\Omega_{m,N}\right)^{p_m}.
\end{equation}
Consequently, no estimate
$\mathfrak B\geq c\,d_{\mathrm{disk}}^{q}$ with $q<p_m$ can hold uniformly on
all centrally symmetric controlled $C^m$ classes.
\end{proposition}

\begin{proof}
Because $N$ is even, $h_{m,N}(\theta+\pi)=h_{m,N}(\theta)$, and the first sine and
cosine coefficients vanish, so the corresponding domains are in Steiner
position. Moreover,
\[
 \rho_{m,N}=h_{m,N}+h_{m,N}''
 =1+\eta(1-N^2)N^{-m}\cos(N\theta).
\]
For fixed sufficiently small $\eta$, the radii $\rho_{m,N}$ stay, for example,
between $1/2$ and $3/2$, while $h_{m,N}>0$ and
$\left\|h_{m,N}\right\|_{C^m}$ is bounded
independently of $N$. The standard support-function parametrization
\eqref{eq:p} therefore produces centrally symmetric strictly convex domains
lying in a single controlled class. For every
constant $r$,
\[
 \left\|h_{m,N}-r\right\|_{C^2}\geq\left\|h_{m,N}''\right\|_\infty
 =\eta N^{-(m-2)},
\]
whereas the choice $r=1$ gives
$\left\|h_{m,N}-1\right\|_{C^2}\leq3\eta N^{-(m-2)}$. This proves the first comparison in
\eqref{eq:sharp-asymptotics}. Since
$\left\|h_{m,N}'\right\|_\infty=\eta N^{-(m-1)}$, the identity
\eqref{eq:symmetric-defect-identity} and the uniform bounds for $\rho_{m,N}$
give the second comparison. The last assertion follows because, for
$q<p_m$,
\[
 \frac{\mathfrak B\left(\Omega_{m,N}\right)}{d_{\mathrm{disk}}\left(\Omega_{m,N}\right)^q}
 \asymp_{m,\eta}N^{-(m-1)+(m-2)q}\longrightarrow0.
\]
\end{proof}
\begin{remark}[What is and is not sharp]\label{rem:sharp-scope}
\emph{Proposition \ref{prop:scale-sharpness} proves sharpness of the geometric
conversion from $C^2$ distance to defect at every finite regularity
$m\geq5$. It does not assert that the spectral cutoff exponents $2p_m$ or
$2$ are attained by actual
overdetermined eigenpairs; the stationary-phase argument gives sufficient
exclusion thresholds. For fixed-mode symmetric perturbations, both
$d_{H,{\mathrm{disk}}}$ and $\mathfrak B$ are first order in the perturbation,
so the linear scale in \eqref{eq:defect-symmetric-lower} is also the correct
small-perturbation scale.}
\end{remark}

In conventional boundary regularity terms, a bounded $C^6$ domain with
strictly positive curvature satisfies the support-function hypothesis. Indeed, a $C^6$ regular boundary parametrization has a $C^5$
Gauss map. Positive curvature makes the lifted normal-angle map a $C^5$
diffeomorphism; its inverse is $C^5$, and therefore the normal-angle
parametrization $p(\theta)$ and
$h(\theta)=p(\theta)\cdot n(\theta)$ are of class $C^5$.

\subsection{Uniform two-point stationary phase}

The normal-angle parametrization \eqref{eq:p} is convenient because
\begin{equation}\label{eq:p-prime}
p'(\varphi)=\rho(\varphi)\tau(\varphi),
\qquad
\dd s=\rho(\varphi)\dd\varphi.
\end{equation}
Fix $\theta$ and set
\[
\xi=n(\theta),\qquad \eta=\tau(\theta).
\]
Define
\begin{equation*}\label{eq:phase}
f_\theta(\varphi)=n(\theta)\cdot p(\varphi).
\end{equation*}
By \eqref{eq:p-prime},
\begin{equation}\label{eq:phase-derivative}
f_\theta'(\varphi)
=\rho(\varphi)\sin(\theta-\varphi).
\end{equation}
Therefore the phase has exactly two critical points,
\[
\varphi=\theta,\qquad \varphi=\theta+\pi.
\]
They are nondegenerate, because
\begin{equation}\label{eq:phase-hessian}
f_\theta''(\theta)=-\rho(\theta),\qquad
f_\theta''(\theta+\pi)=\rho(\theta+\pi).
\end{equation}
Moreover,
\begin{equation}\label{eq:stationary-values}
\begin{aligned}
f_\theta(\theta)&=h(\theta),&
\eta\cdot p(\theta)&=h'(\theta),\\
f_\theta(\theta+\pi)&=-h(\theta+\pi),&
\eta\cdot p(\theta+\pi)&=-h'(\theta+\pi).
\end{aligned}
\end{equation}

We first record the uniform form of the one-dimensional stationary-phase
estimate that will be used below. This formulation makes the dependence of
the remainder constant explicit at the level needed for the geometric
class.

\begin{lemma}[Uniform two-critical-point stationary phase]
\label{lem:uniform-SP}
Let $K,\mu,r>0$. Consider a family of real-valued phases
$f\in C^4\left(\mathbb R/2\pi\mathbb Z\right)$ and complex-valued amplitudes
$a\in C^3\left(\mathbb R/2\pi\mathbb Z\right)$ with exactly two critical points
$q_+$ and $q_-$. Assume that
\begin{equation*}\label{eq:abstract-SP-bounds}
 \|f\|_{C^4}+\|a\|_{C^3}\leq K,
 \qquad \left|f''\right|\geq \mu
 \quad\hbox{on the $2r$-neighborhoods of }q_\pm,
\end{equation*}
and
\begin{equation}\label{eq:abstract-nonstationary}
 \left|f'\right|\geq \mu
 \quad\hbox{outside the $r$-neighborhoods of }q_\pm,
\end{equation}
where the two $2r$-neighborhoods are disjoint. Then there is a constant
$C_{\mathrm{sp}}=C_{\mathrm{sp}}(K,\mu,r)$ such that, for every $\lambda\geq1$,
\begin{align*}
\int_0^{2\pi}a(\varphi)e^{i\lambda f(\varphi)}\,d\varphi
={}&\left(\frac{2\pi}{\lambda}\right)^{1/2}
\sum_{\sigma\in\{+,-\}}
\frac{a(q_\sigma)}{\sqrt{\left|f''(q_\sigma)\right|}}
e^{i\lambda f\left(q_\sigma\right)
+i\frac{\pi}{4}\operatorname{sgn}f''\left(q_\sigma\right)}
+\mathcal E_\lambda,
\label{eq:abstract-SP}
\end{align*}
with
\begin{equation}\label{eq:abstract-SP-error}
 \left|\mathcal E_\lambda\right|\leq C_{\mathrm{sp}}\lambda^{-3/2}.
\end{equation}
\end{lemma}

\begin{proof}
We give the normalization needed to apply the cited finite-regularity result.
Set
\[
 r_0:=\min\left\{\frac r2,\frac{\mu}{8K},\frac18\right\}>0
\]
and choose cutoffs $\chi_\pm$ supported in the $2r_0$-neighborhoods of
$q_\pm$ and equal to one on the corresponding $r_0$-neighborhoods. Their
$C^3$ norms are bounded in terms of $K,\mu,$ and $r$. Fix one critical
point $q$ and put $\sigma=\operatorname{sgn}f''(q)$. For $|s|\le2r_0$,
\[
 \left|f''(q+s)-f''(q)\right|\le 2Kr_0\le\frac\mu4
 \le\frac14\left|f''(q)\right|.
\]
Consequently $f''$ has the fixed sign $\sigma$ there. With
\[
 \gamma_q=\left(2r_0\right)^2\left|f''(q)\right|,
 \qquad
 \Phi_q(x)=
 \frac{\sigma\bigl(f\left(q+2r_0x\right)-f(q)\bigr)}{\gamma_q},
 \qquad |x|\le1,
\]
we have
\[
 \Phi_q(0)=\Phi_q'(0)=0,
 \qquad
 \frac34\le\Phi_q''(x)\le\frac54.
\]
Moreover, $\left\|\Phi_q\right\|_{C^4}$ is bounded by a constant depending only on
$K,\mu,$ and $r$. After the change of variables
$\varphi=q+2r_0x$, the localized integral is
\[
 e^{i\lambda f(q)}
 \int_{-1}^{1}\Psi_q(x)
 e^{i\sigma\Lambda_q\Phi_q(x)}\dd x,
 \quad
 \Psi_q(x)=2r_0\chi_q\left(q+2r_0x\right)a\left(q+2r_0x\right),
 \quad
 \Lambda_q=\lambda\gamma_q.
\]
The $C^2$ norm of $\Psi_q$ is bounded by $C(K,\mu,r)$, and
\begin{equation}\label{eq:effective-parameter}
 \mu\left(2r_0\right)^2\lambda
 \le\Lambda_q\le K\left(2r_0\right)^2\lambda.
\end{equation}

The one-dimensional estimate \cite[Theorem 2.1(2)]{Guo} applies directly to
this normalized integral when $\sigma=1$. When $\sigma=-1$, we use the exact
identity
\[
 \int_{-1}^{1}\Psi_q(x)e^{-i\Lambda_q\Phi_q(x)}\dd x
 =\overline{\int_{-1}^{1}\overline{\Psi_q(x)}
 e^{i\Lambda_q\Phi_q(x)}\dd x}
\]
and apply the same estimate to the integral on the right-hand side before
taking complex conjugates. Thus the remainder bound is unchanged and the
stationary-phase factor is $e^{-i\pi/4}$.
It requires only the derivatives of the phase through order four and of the
amplitude through order two, and its remainder constant is controlled by
those norms under the displayed uniform lower and upper bounds for
$\Phi_q''$. Using \eqref{eq:effective-parameter} and
$\Psi_q(0)=2r_0a(q)$ gives
\begin{align*}
 \int \chi_q(\varphi)a(\varphi)e^{i\lambda f(\varphi)}\dd\varphi
 ={}&\left(\frac{2\pi}{\lambda}\right)^{1/2}
 \frac{a(q)}{\sqrt{|f''(q)|}}
 e^{i\lambda f(q)+i\pi\sigma/4}
 +O_{K,\mu,r}(\lambda^{-3/2}).
\end{align*}
When $\Lambda_q<1$, \eqref{eq:effective-parameter} confines $\lambda$ to a
fixed compact interval, and the same estimate follows after increasing the
constant. Although the cited theorem is stated for smooth data, the stated
$C^4$ phase and $C^3$ amplitude classes follow by approximation in
$C^4\times C^2$: after mollification one subtracts the affine Taylor part of
the phase at the origin, preserving the normalized second-derivative bounds,
and then passes to the limit. Thus no unlisted higher smoothness norm enters
the constant. It is worth noting that Proposition~\ref{Proposition 2.1} gives the same leading coefficient.

It remains to estimate the complement of the two localized pieces. Inside
the original $r$-neighborhoods but outside the $r_0$-neighborhoods, the fixed
sign of $f''$ and $\left|f''\right|\ge\mu$ imply $\left|f'\right|\ge\mu r_0$; outside the original
$r$-neighborhoods, \eqref{eq:abstract-nonstationary} applies. Hence $\left|f'\right|$
has a positive lower bound depending only on $\mu$ and $r_0$ throughout the
complement. Two integrations by parts with
\[
 \frac{1}{i\lambda f'(\varphi)}\frac{\mathrm{d}}{\mathrm{d}\varphi}
\]
give $O\left(\lambda^{-2}\right)$; the resulting coefficients are bounded in terms of
$K$, $\mu$, and $r$. Adding the two localized expansions and increasing one
constant gives
\eqref{eq:abstract-SP-error}.
\end{proof}

Based on Lemma \ref{lem:uniform-SP}, we can establish the following two-point expansion.
\begin{proposition}[Uniform two-point expansion]\label{prop:stationary}
Let $\Omega\in\cG_5\left(R,M,\kappa_0,\kappa_1\right)$. There exists
\[
C_*=C_*\left(R,M,\kappa_0,\kappa_1\right)>0
\]
such that, uniformly for $\theta\in\R$, $|t|\le1$, and $\lambda\ge3$,
\begin{align}
I_D(\lambda,t,\theta)
&:=\int_{\pd\Omega}
e^{t\eta\cdot x}e^{i\lambda\xi\cdot x}\dd s\nonumber\\
&=\left(\frac{2\pi}{\lambda}\right)^{1/2}
\left[
\sqrt{\rho(\theta)}e^{t h'(\theta)}
e^{i\lambda h(\theta)-i\pi/4}
+\sqrt{\rho(\theta+\pi)}e^{-t h'(\theta+\pi)}
e^{-i\lambda h(\theta+\pi)+i\pi/4}
\right]
+\cR_D,\label{eq:SP-D}\\
I_N(\lambda,t,\theta)
&:=\frac1{i\lambda}
\int_{\pd\Omega}
\pd_\nu\!\left(
e^{t\eta\cdot x}e^{i\lambda\xi\cdot x}
\right)\dd s\nonumber\\
&=\left(\frac{2\pi}{\lambda}\right)^{1/2}
\left[
\sqrt{\rho(\theta)}e^{t h'(\theta)}
e^{i\lambda h(\theta)-i\pi/4}
-\sqrt{\rho(\theta+\pi)}e^{-t h'(\theta+\pi)}
e^{-i\lambda h(\theta+\pi)+i\pi/4}
\right]
+\cR_N.\label{eq:SP-N}
\end{align}
The remainders satisfy
\begin{equation}\label{eq:remainder}
\left|\cR_D\right|+\left|\cR_N\right|\le C_*\lambda^{-3/2}.
\end{equation}
\end{proposition}

\begin{proof}
By \eqref{eq:p-prime},
\begin{equation*}\label{eq:I-D-param}
I_D
=\int_0^{2\pi}
a_{D,t,\theta}(\varphi)
e^{i\lambda f_\theta(\varphi)}\dd\varphi,
\quad
a_{D,t,\theta}(\varphi)
=\rho(\varphi)e^{t\tau(\theta)\cdot p(\varphi)}.
\end{equation*}
For the Neumann integral,
\begin{align*}
\frac1{i\lambda}\pd_\nu
\left(e^{t\eta\cdot x}e^{i\lambda\xi\cdot x}\right)
&=\left(
\xi\cdot\nu-i\frac{t}{\lambda}\eta\cdot\nu
\right)
e^{t\eta\cdot x}e^{i\lambda\xi\cdot x}.
\label{eq:N-amplitude}
\end{align*}
Since $\nu(p(\varphi))=n(\varphi)$, its amplitude, including arclength, is
\begin{equation}\label{eq:aN}
a_{N,t,\theta,\lambda}(\varphi)
=\rho(\varphi)
\left(
n(\theta)\cdot n(\varphi)
-i\frac{t}{\lambda}\tau(\theta)\cdot n(\varphi)
\right)
e^{t\tau(\theta)\cdot p(\varphi)}.
\end{equation}

We make the finite-regularity dependence explicit. From
\eqref{eq:p} and $\rho=h+h''$,
\begin{equation*}\label{eq:p-rho-regularity}
 \|p\|_{C^4}\le C\|h\|_{C^5}\le CM,
 \qquad
 \|\rho\|_{C^3}\le2M,
 \qquad
 \|p\|_{L^\infty}\le R.
\end{equation*}
Consequently, uniformly for $\theta\in\R$ and $|t|\le1$,
\begin{equation*}\label{eq:exponential-amplitude-bound}
 \left\|e^{t\tau(\theta)\cdot p(\cdot)}\right\|_{C^3}
 \le C e^R(1+M)^3.
\end{equation*}
The fixed trigonometric factors in \eqref{eq:aN}, together with
$|t|/\lambda\le1$, then give
\[
 \left\|f_\theta\right\|_{C^4}\le CM,
 \qquad
 \left\|a_{D,t,\theta}\right\|_{C^3}
 +\left\|a_{N,t,\theta,\lambda}\right\|_{C^3}
 \le C\left(R,M,\kappa_0,\kappa_1\right).
\]
Thus
\begin{equation*}\label{eq:phase-amplitude-K}
\sup_{\substack{\theta\in\R,\ |t|\leq1\\ \lambda\geq1}}
\left(
\left\|f_\theta\right\|_{C^4}
+\left\|a_{D,t,\theta}\right\|_{C^3}
+\left\|a_{N,t,\theta,\lambda}\right\|_{C^3}
\right)
\leq K_*=K_*\left(R,M,\kappa_0,\kappa_1\right).
\end{equation*}

Since $\rho'=h'+h'''$, we have $\left\|\rho'\right\|_{L^\infty}\leq2M$.  Choose
$r_*=r_*\left(M,\kappa_1)\in(0,\pi/8\right)$ so that
\begin{equation*}\label{eq:r-star}
 \cos\left(2r_*\right)\geq\frac34,
 \qquad 2M\sin\left(2r_*\right)\leq\frac{1}{4\kappa_1}.
\end{equation*}
Writing $s=\varphi-\theta$ near $\theta$, and
$s=\varphi-\theta-\pi$ near $\theta+\pi$, the identity
\eqref{eq:phase-derivative} gives
\begin{equation*}\label{eq:local-hessian-uniform}
 \left|f_\theta''(\varphi)\right|
 \geq \frac{1}{\kappa_1}\cos\left(2r_*\right)
       -2M\sin\left(2r_*\right)
 \geq\frac{1}{2\kappa_1}
\end{equation*}
on both $2r_*$-neighborhoods. Outside the corresponding
$r_*$-neighborhoods,
\begin{equation*}\label{eq:nonstationary-uniform}
 \left|f_\theta'(\varphi)\right|
 =\rho(\varphi)|\sin(\theta-\varphi)|
 \geq \frac{\sin r_*}{\kappa_1}.
\end{equation*}
Thus Lemma \ref{lem:uniform-SP} applies with
\begin{equation*}\label{eq:mu-star}
 \mu_*=\min\left\{\frac{1}{2\kappa_1},
                  \frac{\sin r_*}{\kappa_1}\right\}.
\end{equation*}

For the Dirichlet amplitude, \eqref{eq:stationary-values} gives
\[
\frac{a_{D,t,\theta}(\theta)}{\sqrt{\rho(\theta)}}
=\sqrt{\rho(\theta)}e^{t h'(\theta)}
\]
and
\[
\frac{a_{D,t,\theta}(\theta+\pi)}
{\sqrt{\rho(\theta+\pi)}}
=\sqrt{\rho(\theta+\pi)}e^{-t h'(\theta+\pi)}.
\]
Together with \eqref{eq:phase-hessian}, Lemma \ref{lem:uniform-SP} proves
\eqref{eq:SP-D}, with remainder constant bounded by
$C_{\mathrm{sp}}\left(K_*,\mu_*,r_*\right)$.

For the Neumann amplitude \eqref{eq:aN},
\[
n(\theta)\cdot n(\theta)=1,\qquad
n(\theta)\cdot n(\theta+\pi)=-1,
\]
while the term multiplied by $t/\lambda$ vanishes. A second application of
Lemma \ref{lem:uniform-SP} gives \eqref{eq:SP-N}. Taking, for example,
\[
 C_*=2C_{\mathrm{sp}}\left(K_*,\mu_*,r_*\right)
\]
gives \eqref{eq:remainder} and shows that $C_*$ depends only on
$R,M,\kappa_0,$ and $\kappa_1$, but not on $\delta$, $\Omega$, $\theta$,
$t$, or $\lambda$.
\end{proof}

\begin{remark}[Signed curvature]
\emph{Uniform convexity gives $\kappa>0$. The leading factor is therefore
$\kappa^{-1/2}=\sqrt{\rho}$, while the signs of the phase Hessians are
recorded separately by $e^{-i\pi/4}$ and $e^{i\pi/4}$. This avoids an
ambiguous square root of signed curvature.}
\end{remark}

\subsection{The Dirichlet problem with constant Neumann trace}

\begin{lemma}[Dirichlet high-frequency exclusion]\label{thm:D}
Under the assumptions of Proposition \ref{thm:uniform}, problem \eqref{overdeterminedeigenvalueproblem}
has no nontrivial solution whenever
\begin{equation}\label{eq:threshold-D}
\alpha>
\left[
\max\left\{
3,\,
\frac{C_*}{\sqrt{2\pi}\,\delta}
\right\}
\right]^2-1.
\end{equation}
\end{lemma}

\begin{proof}
By \eqref{eq:defect}, choose $\theta_0$ and $t_0\in\{-1,1\}$ such that
\[
\left|
\sqrt{\rho\left(\theta_0\right)}e^{t_0h'\left(\theta_0\right)}
-\sqrt{\rho\left(\theta_0+\pi\right)}e^{-t_0h'\left(\theta_0+\pi\right)}
\right|\geq\delta.
\]
Set $\lambda=\sqrt{\alpha+1}$. Then \eqref{eq:dispersion} holds, and
Lemma \ref{prop:green} gives
$I_D\left(\lambda,t_0,\theta_0\right)=0$. Hence
\eqref{eq:SP-D} and the reverse triangle inequality give
\[
\left|\cR_D\right|
\ge
\left(\frac{2\pi}{\lambda}\right)^{1/2}
\left|
\sqrt{\rho\left(\theta_0\right)}e^{t_0h'\left(\theta_0\right)}
-\sqrt{\rho\left(\theta_0+\pi\right)}e^{-t_0h'\left(\theta_0+\pi\right)}
\right|
\geq\left(\frac{2\pi}{\lambda}\right)^{1/2}\delta.
\]
On the other hand, \eqref{eq:remainder} yields
$\left|\cR_D\right|\le C_*\lambda^{-3/2}$. Therefore
\begin{equation}\label{eq:lambda-upper-D}
\lambda
\le
\frac{C_*}{\sqrt{2\pi}\,\delta}.
\end{equation}
Condition \eqref{eq:threshold-D} contradicts
\eqref{eq:lambda-upper-D}.
\end{proof}

\begin{remark}
\emph{No phase alignment is used. The complex phases of the two leading terms may
be arbitrary; the reverse triangle inequality depends only on their unequal
moduli. The defect \eqref{eq:defect} measures this inequality directly and
does not require the two curvature factors to be equal.}
\end{remark}

\begin{proof}[Proof of Proposition~\ref{thm:uniform}]
Take
\begin{equation}\label{eq:A-uniform}
A\left(R,M,\kappa_0,\kappa_1,\delta\right)
=
\left[
\max\left\{
3,\,
\frac{C_*\left(R,M,\kappa_0,\kappa_1\right)}
{\sqrt{2\pi}\,\delta}
\right\}
\right]^2-1.
\end{equation}
Lemma~\ref{thm:D} proves the result.
\end{proof}

Now we present the proof of Theorem 1.1.

\begin{proof}[Proof of Theorem 1.1]
We would like to mention that $\partial \Omega$ is real analytic if $\partial \Omega$ is $C^{2,\epsilon}$.
Given that $\partial \Omega$ is of class $C^{2,\epsilon}$, we can deduce from \cite[Theorem 6.14]{Gilbarg} that $u\in C^2\left(\overline{\Omega}\right)$.
Additionally, the regularity result \cite[Theorem 2]{Kinderlehrer} remains valid if the origin is replaced by any fixed point on $\partial\Omega$.
We fix a point on $\partial\Omega$ where the unit normal vector points inward as $(0,\pm1)$.
Consider $g(p)=|p|^2-c^2$. At this fixed point, we have $\left\vert \nabla u\right\vert \neq0$ and $\frac{\partial g}{\partial p_n}\left(\nabla u\right)\neq0$.
Applying \cite[Theorem 2]{Kinderlehrer}, we conclude that the boundary $\partial \Omega$ is real analytic. So, $\Omega$ automatically belongs to the set class defined in Definition~\ref{def:class}.

Translate $\Omega$ into Steiner position; this changes neither
overdetermined problem. Since $\Omega$ is not a disk,
\[
\beta_\Omega:=\mathfrak B(\Omega)>0
\]
by Proposition \ref{prop:defect-disk}.
Choose $R_\Omega$, $M_\Omega$, and positive curvature bounds
$\kappa_{0,\Omega}$ and $\kappa_{1,\Omega}$ for this fixed domain.
Proposition \ref{thm:uniform} with $m=5$ and
$\delta=\beta_\Omega/2$ defines
\begin{equation*}\label{eq:alpha0}
\alpha_0(\Omega)
=A\left(
R_\Omega,M_\Omega,
\kappa_{0,\Omega},\kappa_{1,\Omega},
\frac12\beta_\Omega
\right).
\end{equation*}
An overdetermined eigenvalue above this number would contradict Proposition
\ref{thm:uniform}, proving the asserted exclusion.
\end{proof}

\begin{remark}[Scope of the conclusion]\label{rem:scope}
\emph{The threshold in Theorem 1.1 depends on the full geometry of
the fixed domain. It is not a universal threshold over all smooth planar
domains and does not contradict the counterexamples in
\cite{ColbrookStepaniants,CaoDeDios}. Moreover, the qualitative existence of
a fixed-domain cutoff is consistent with the classical sequence-rigidity
results: a non-disk cannot carry an unbounded sequence of distinct
overdetermined eigenvalues. Proposition \ref{thm:main-distance} gives the
nonsymmetric $C^2$ statement with regularity-dependent cutoff
$C_m(1+\varepsilon^{-2p_m})$, whereas Proposition \ref{thm:symmetric-H} gives the
centrally symmetric Hausdorff statement with cutoff
$C(1+\varepsilon^{-2})$. Both are high-frequency results; neither excludes
the finitely many eigenvalues below its threshold. The defect-based
Proposition \ref{thm:uniform} is the common analytic input.}
\end{remark}

\begin{remark}[Comparison with sequence results]
\emph{The classical argument of Berenstein and Yang \cite{Berenstein2} uses an
unbounded sequence of overdetermined eigenvalues together with a slowly
growing exponential parameter. Our argument uses one eigenvalue and
uniform convexity.  The defect \eqref{eq:defect} replaces both a frequency
limit and the former central-symmetry assumption. This does not turn the
fixed-domain theorem into an all-frequency rigidity theorem: finitely many
exceptional eigenvalues below the domain-dependent cutoff are not excluded.}
\end{remark}

\section{The proof of Theorem 1.2}
\bigskip
\quad\,
Due to the similarity between the proof of Theorem 1.2 and Theorem 1.1, we only outline the general approach and highlight the key differences.
The boundary
integral for the Neumann problem contains the variable normal factor
$(t\eta+i\lambda\xi)\cdot\nu$. The normalization by $i\lambda$ in
\eqref{eq:SP-N} keeps this amplitude uniform and records its values $1$ and
$-1$ at the two stationary points.
\begin{lemma}[Schiffer high-frequency exclusion]\label{thm:N}
Under the assumptions of Proposition \ref{thm:uniform}, problem \eqref{SchifferConjecture}
has no nontrivial solution whenever
\begin{equation}\label{eq:threshold-N}
\alpha>
\left[
\max\left\{
3,\,
\frac{C_*}{\sqrt{2\pi}\,\delta}
\right\}
\right]^2-1.
\end{equation}
\end{lemma}

\begin{proof}
Choose $\theta_0$, $t_0$, and $\lambda$ exactly as in Lemma \ref{thm:D}.
Lemma \ref{prop:green} gives
$I_N(\lambda,t_0,\theta_0)=0$. The additional minus sign in
\eqref{eq:SP-N} changes only the relative complex phase, not the difference
of the two moduli. Hence
\[
|\cR_N|
\ge
\left(\frac{2\pi}{\lambda}\right)^{1/2}
\delta.
\]
Combining this with $|\cR_N|\le C_*\lambda^{-3/2}$ gives
\[
\lambda
\le
\frac{C_*}{\sqrt{2\pi}\,\delta},
\]
contradicting \eqref{eq:threshold-N}.
\end{proof}

Since $\partial \Omega$ is Lipschitz, by \cite{Caffarelli, Williams0, Williams}, the boundary $\partial \Omega$ is real analytic. Thus, $\Omega$ belongs to the set class defined in Definition~\ref{def:class}. Then following the same steps as that of Theorem 1.1 with the appropriate modifications, we can complete the proof of Theorem 1.2.

\section{Sharpness and scope}
\bigskip
\quad\,
The proof separates into an analytic module and two geometric conversions.
The analytic module, Proposition \ref{thm:uniform}, uses only a controlled
$C^5$ support function. Uniform convexity gives exactly two nondegenerate
stationary points, while the radius, curvature, and $C^5$ bounds make the
$O(\lambda^{-3/2})$ remainder uniform. The factor $\delta^{-1}$ in
\eqref{eq:A-uniform} records the possible approach of the two leading
moduli. Since the argument fixes $t=\pm1$, the dispersion relation is
$\alpha=\lambda^2-1$; this explains why a defect lower bound of size
$\delta$ produces an eigenvalue cutoff of size $O\left(1+\delta^{-2}\right)$.

Without symmetry, converting geometric distance into defect requires the
even--odd decomposition and the $C^m$ interpolation estimate
\eqref{eq:even-C2-interpolation}. At finite regularity $m\geq5$, this gives
the defect exponent $p_m=(m-1)/(m-2)$ and the eigenvalue cutoff exponent
$2p_m$. These exponents decrease to $1$ and $2$, respectively, as
$m\to\infty$. Proposition \ref{prop:scale-sharpness} shows that $p_m$ cannot
be lowered under only uniform $C^m$ control, even within the centrally
symmetric subclass.

Its real benefit is instead the exact identity
\eqref{eq:symmetric-defect-identity}.  It removes the interpolation step and
converts the weaker centered Hausdorff distance linearly into defect already
under $C^5$ control. The resulting cutoff has order
$1+\varepsilon^{-2}$. Keeping the two metrics separate is essential: the
regularity-dependent theorem is stronger in geometric scope, while the
symmetric theorem is stronger in metric and cutoff exponent.

None of the constants is universal. If a curvature bound degenerates, the
controlled norm diverges, or the relevant distance from disks tends to
zero, the threshold may diverge. The counterexamples in
\cite{ColbrookStepaniants,CaoDeDios} therefore do not conflict with the
results: the present theorems exclude only sufficiently high frequencies
after a controlled class and a positive distance from disks have been
fixed. Conversely, the conclusions are not all-frequency versions of the
Schiffer conjecture and leave open the eigenvalues below the stated cutoff.
\\ \\
\textbf{The conflicts of interest statement and Data Availability statement.}
\bigskip\\
\indent There is not any conflict of interest.
Data sharing not applicable to this article as no datasets were generated or analysed during the current study.
\\ \\
\textbf{AI declaration.}
\bigskip\\
\indent ChatGPT 5.6 in the form of codex was used to help debug an early version
of some of the code. All of its outputs were checked line by line by the authors. The authors
take full responsibility for all content.
\bibliographystyle{amsplain}
\makeatletter
\def\@biblabel#1{#1.~}
\makeatother


\providecommand{\bysame}{\leavevmode\hbox to3em{\hrulefill}\thinspace}
\providecommand{\MR}{\relax\ifhmode\unskip\space\fi MR }
\providecommand{\MRhref}[2]{%
  \href{http://www.ams.org/mathscinet-getitem?mr=#1}{#2}
}
\providecommand{\href}[2]{#2}

\end{document}